# Spatial-temporal data mining procedure: LASR

## Xiaofeng Wang[1], Jiayang Sun[2] and Kath Bogie[3]

*The Cleveland Clinic Foundation, Case Western Reserve University and Cleveland FES center*

**Abstract:** This paper is concerned with the statistical development of our spatial-temporal data mining procedure, LASR (pronounced "laser"). LASR is the abbreviation for *Longitudinal Analysis with Self-Registration* of large-$p$-small-$n$ data. It was motivated by a study of "Neuromuscular Electrical Stimulation" experiments, where the data are noisy and heterogeneous, might not align from one session to another, and involve a large number of multiple comparisons. The three main components of LASR are: (1) data segmentation for separating heterogeneous data and for distinguishing outliers, (2) automatic approaches for spatial and temporal data registration, and (3) statistical smoothing mapping for identifying "activated" regions based on false-discovery-rate controlled $p$-maps and movies. Each of the components is of interest in its own right. As a statistical ensemble, the idea of LASR is applicable to other types of spatial-temporal data sets beyond those from the NMES experiments.

## 1. Introduction

Developments of medical and computer technology in the last two decades have enabled us to collect huge amounts of data in both spatial and temporal dimensions. These types of data have become common in medical imaging, epidemiology, neuroscience, ecology, climatology, environmentology and other areas. Typical spatial-temporal data can be denoted by $y(s, t, i)$, where $y$ is the intensity value at the spatial location $s \in \mathcal{S}$, time $t \in \mathcal{T}$ and for the subject indexed by $i \in \mathcal{N}$. In most applications, $\mathcal{S} = \{1, \ldots, S\}$ is a 1, 2 or 3 dimensional space indexed by $S$ pixels; $\mathcal{T} = \{1, 2, \ldots, T\}$ is a set of $T$ time points; and $\mathcal{N} = \{1, 2, \ldots, n\}$ is the set of $n$ subjects. In principle, the indexing can be done by continuous variables, but in practice, only a discretized version is observed. It is often the case that the data size $n$ is much smaller than the data dimension $p = S \times T$. Hence the data are of *large-p-small-n*.

An example of such spatial-temporal data is the data from our *Neuromuscular Electrical Stimulation* (NMES) experiments to prevent pressure sores. *Pressure sores* (also called pressure ulcers, bed sores, or decubitus ulcers) [3] are areas of injured skin and tissue. They are usually caused by sitting or lying in one position for long periods of time. This puts pressure on certain areas of the body which in turn reduce the blood supply to the skin and the tissue under the skin and hence a sore may form. Pressure sores are known to be a multi-factor complication that occurs in many wheelchair users due to reduced mobility, e.g., those with spinal cord injury (SCI).

[1]Department of Quantitative Health Sciences, The Cleveland Clinic Foundation, 9500 Euclid Avenue / Wb4, Cleveland, OH 44195, e-mail: wangx6@ccf.org

[2]Department of Statistics, Case Western Reserve University, 10900 Euclid Avenue, Cleveland, Ohio 44106, e-mail: jsun@case.edu

[3]Cleveland FES Center, Hamann Building, Room 601, 2500 MetroHealth Drive, Cleveland, OH 44109-1998; Department of Orthopedics, Case Western Reserve University, 10900 Euclid Avenue, Cleveland, Ohio 44106, e-mail: kath.bogie@case.edu

Traditionally, techniques to reduce pressure sore incidence have focused on extrinsic risk factors by providing cushions which improve pressure distribution and educating individuals on the importance of regular pressure relief procedures. There remains a significant number of people with SCI for whom pressure relief cushions are inadequate and/or who are unable to maintain an adequate pressure relief regime. NMES provides a unique technique to produce beneficial changes at the user/support system interface by altering the intrinsic characteristics of the user's paralyzed tissue itself [4]. To quantify the effects of long-term NMES on the intrinsic characteristics of the paralyzed muscles, data on the response to loading, including interface pressure distribution when seated in a wheelchair, must be acquired over a long period of time and be statistically analyzed.

In Section 2, we describe the background and challenges of data analysis from our NMES experiments that motivated us to develop LASR. In Section 3, we address the important data preprocessing issues in data mining. Two steps are proposed here, data segmentation and data registration. An optimal threshold method with the EM algorithm is proposed to classify the sitting (signal) region from the background in data frames. We then introduce a self-registration technique, *Self-Registration by a Line and a Point* (SRLP) for spatial registration, incorporated by a fast temporal registration scheme, *Intensity-based Correlation Registration* (ICR). In Section 4, we propose a *Statistical Smoothing Mapping* (SSM) algorithm for interface pressure analysis, which includes the multivariate smoothing techniques. Since the number of significance tests for testing the difference regions is equal to the number of pixels per frame, an overall error rate of the tests must be controlled. Here we choose to develop false-discovery-rate (FDR) controlled movies and maps, called FDR movies and FDR maps, under dependency, to overcome the multiplicity effect from testing "activation" pixels simultaneously. In Section 5, combining the techniques in the previous sections, we present a data-mining scheme, the LASR procedure for analyzing a large sequence of spatial-temporal data sets. LASR is shown to be effective in the application to data from the NMES experiments. In Section 6, a discussion on applications of LASR to other fields and future research is given.

## 2. Experimental data and challenges

*Background.* The primary hypothesis of our clinical study is that chronic use of NMES improves pressure distribution at the seating support area, specifically by the reduction of peak pressures over bony prominences. In addition, chronic NMES will increase vascularity leading to improved tissue blood flow and resulting in improved regional tissue health in individuals with SCI. Therefore, repeated assessments of sitting interface pressures were obtained for a group of eight subjects with SCI participating in a study to investigate the use of NMES for standing and transfers. All subjects were full-time wheelchair users at entrance into the study and had sustained traumatic SCI from 13-204 months prior to enrollment. All subjects had complete SCI and were therefore considered to be at increased risk of tissue breakdown, in part due to disuse muscle atrophy of the glutei.

Seating interface pressures were determined using a Tekscan Advanced Clinseat Pressure Mapping System (Tekscan Inc., Boston, Massachusetts). Assessments were carried out prior to commencing regular use of stimulation, to obtain a baseline value, and then at intervals of 3–12 months during their participation in the study, giving an overall time frame of up to five years for repeated assessments of each study participant.

In order to perform an assessment of seating interface pressures the subject transferred out of the wheelchair and a pressure sensor mat was placed over the wheelchair cushion. The sensor mat is comprised of a matrix of pressure sensitive cells (38 rows, 41 columns). The subject then transferred back into the wheelchair and was asked to sit in their customary sitting posture. Care was taken to insure that the sensor mat was not creased or folded under the subject in order to avoid inaccurate high spots. The sensor was then calibrated based on the assumption that 80 percent of the subject's body weight was acting through the seat base. Calibration took less than 20 seconds to complete. Interface pressure data was then collected for 200 seconds at a rate of 2 frames/sec. The subject was then asked to perform a pressure relief procedure and sit back in the same position. The sensor was then recalibrated and a second set of pressure data was collected at the same rate of data collection while left/right alternating gluteal stimulation was applied to provide dynamic side-to-side weight shifting for 200 seconds. Interface pressure data was collected concurrently at a rate of 2 frames/sec. Stimulation was then discontinued and subjects were asked to repeat the pressure relief procedure and sit back in the same position before collecting a third set of interface pressure data with subjects in a quiet sitting posture. Real-time two-dimensional pressure intensity data at the seating interface were produced with the use of the assessment device.

*Data.* In summary, for each subject in each session done over time our data sets consist of *three sub-data sets* each of which is under one of three subsequent assessment conditions: no stimulation, on-off alternation stimulation, and no stimulation, as shown in Figure 1. Each of the sub-data sets consists of a sequence of *400 data frames*. Each data frame represents spatial pressure intensity over the *sitting interface (with S= 38 × 41 pixels)* at a certain time point. Hence the data size is $n = 8$ and the data dimension is at least $p = S \times T = (38 \times 41) \times (400 \times 3 \times 3) = 5,608,800$ where the last number 3 is the minimal number of the sessions we had for the patients who received treatments.

FIG 1. *Data structure in the NMES experiment. There are three sub-data sets in each of assessment sessions, under condition: no stimulation, on-off alternation stimulation, and no stimulation, respectively. Each of the sub-data sets consists of 400 data frames, and may be called a segment.*

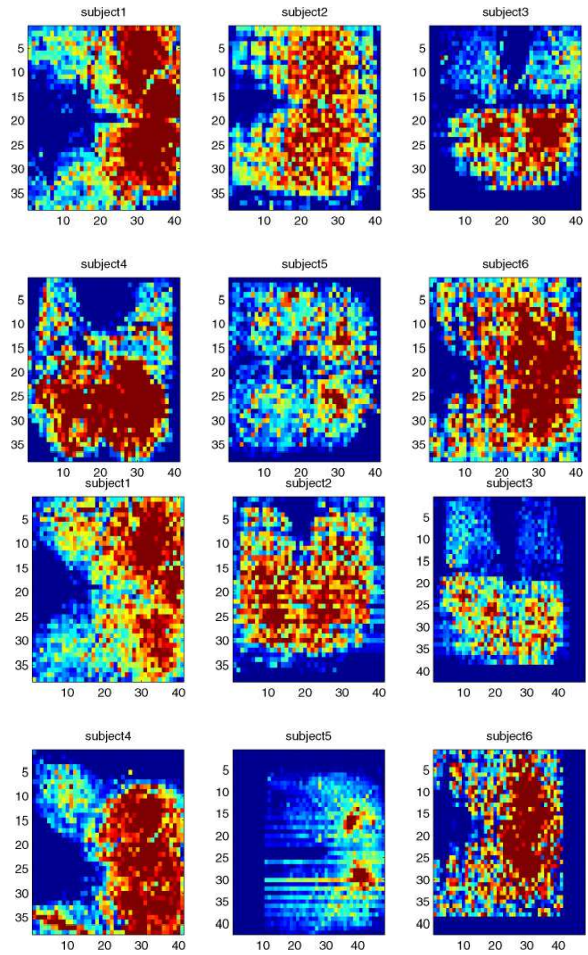

FIG 2. *Top 6 panels are raw data images for the subjects before treatment; and bottom 6 panels are corresponding images after treatment.*

Figure 2 displays the first data frame from each of six sub-data sets (representing six subjects) at the first segment of the first session (before treatment) and at the third segment of the last session (after treatment). Two other subjects are used as control subjects. The numbers of columns and rows correspond to spatial coordinates of a subject's sitting interface. In the movie representation, the $x$-axis and $y$-axis in three-dimensional Cartesian coordinate system denote the spatial coordinates of the sitting interface of subjects; the $z$-axis denotes the pressure intensities. A movie can be generated for each sub-data set. One can then easily see the dynamic changes of pressure intensities. Examples of movies can be found at http://stat.case.edu/lasr, in MPEG format.

The left picture of Figure 3 shows idealized changes in pressure contour across the region of the ischial tuberosities. This is based on comparison *with no electrical stimulation*. Note that the baseline contour shows high mean interface pressures bilaterally in the ischial region which indicates a high risk of local tissue breakdown.

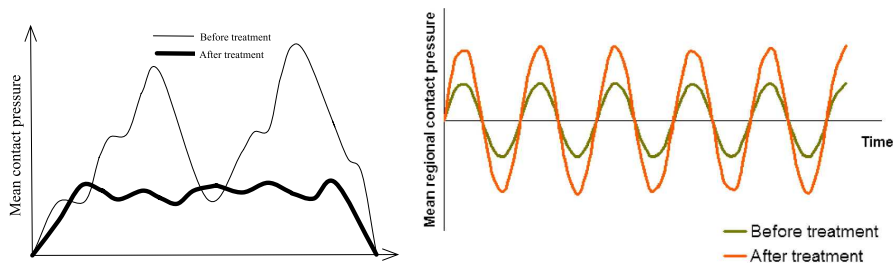

Fig 3. *Left: Idealized changes in pressure contour across the region of the ischial tuberosities under no stimulation; Note the difference in the pressure distributions. Right: Idealized changes in gluteal pressure variation with electrical stimulation over time. Note the change in amplitudes.*

Improved pressure distribution with reduced ischial region interface pressures and more evenly distributed seating pressures indicates a lower risk of tissue breakdown.

Clinicians are also interested in exploring the dynamic changes of interface pressure distribution during electrical stimulation. The right picture of Figure 3 displays idealized changes in pressure variation based on comparison *with electrical stimulation over time*. Regional interface pressures vary cyclically with applied stimulation before treatment. Variations about the mean increase in amplitude (after treatment) because of increasing strength of muscle contraction after long-term treatment. In order to show whether this objective has been met over time and/or with different seating setups there must be some basis for comparison between measurements, so that true differences can be determined.

*Challenges.* (1) *Segmentation and registration for a large sequence of data frames.* In a data mining process, raw data often require some initial processing in order to become useful for further statistical inference, e.g. filtering, scaling, calibration etc. In our NMES study: i) Raw data frames contain background noises; ii) Data frames recorded at different sessions over time from the same subject may not align spatially because, either the subject did not sit in the same relative position on the sensor mat or with the same posture at each assessment, or the image target regions differ from one session to another; iii) Artificial differences between alternating left/right simulation responses can obscure true differences if the data frames from different phases of the stimulation cycle are not aligned temporally between sessions.

Registration techniques have been well developed in the medical imaging area [7, 10]. However, most existing image registration procedures require a reference image and a similarity measure for each candidate image. They are not efficient for calibrating a large number of spatial-temporal data sets, such as registering sequences of data frames or movies in pressure mapping. It would be "labor intensive" to identify the landmarks one by one for each data frame if we used corresponding landmark-based registration for thousands of data frames. Developing effective and fast spatial and temporal registration/calibration algorithms for a large volume of spatial-temporal data sets is important. In Section 3, we first develop a *segmentation* procedure and then a spatial and temporal registration procedure.

(2) *Analysis of large-p-small-n data.* The experimental protocol for this NMES study produced many time points and three assessment conditions for each subject. Thus, the data obtained from the NMES experiment exhibit a *large-p-small-n* problem; that is, a large number of features (pressure intensities) over space and time relative to a small number of subject samples. As given in Section 2, $p$ is greater than 5 million and $n = 8$. Traditional statistical approaches usually are based on the assumption that $p < n$ and are not applicable here without a "transformation". Here, we resolve this problem by performing the subject-by-subject comparisons based on the before and after treatment differences. The differences (after registrations) at these pixel values (frame-by-frame) will become a difference movie and will be treated as if they were regression data points. Hence we have literally transferred $S$ pixels into the "subjects" domain, with now $n = S \times 8$ and $p \leq 400$, and hence a statistical smoothing mapping can be developed – See details in Section 4.

## 3. Data preprocessing

We propose two procedures for pre-processing raw data in this section: (1) *Data Segmentation* for data cleaning; and (2) *Data Registration* for data calibration. Segmentation is important here in that it makes the next step, registration based on random landmarks (estimated from data), more robust. Registration is the process of transforming the different sets of data into one coordinate system. Registration is necessary both spatially and temporally in order to be able to compare and model the data obtained at different times and from different perspectives that are in different coordinate systems.

### 3.1. Data segmentation

As shown in Figure 2, noise and outliers appeared outside of the sitting region (*i.e.* the buttock and thigh region). It is critical to detect the edge of the sitting region of subjects by segmenting the data into the spatial regions of interest and the background in each frame and to remove the background noise by zeroing the corresponding values, before automatically estimating the landmarks illustrated in the next registration step.

We propose a density-based segmentation method in which a pixel will be classified into the background, ie the non-sitting region, if its intensity value is less than a threshold $T$. Let $Z(i,j)$ denote the intensity value of the $i$th row and the $j$th column of a data frame. Then the segmented image will have the intensity values:

$$\tilde{Z}(i,j) = \begin{cases} Z(i,j), & \text{if } Z(i,j) > T; \\ 0, & \text{if } Z(i,j) \leq T. \end{cases}$$

A simple and effective way of computing $T$ is to model the density of intensity values at all pixels in each frame by a mixture of normal distributions:

$$(1) \qquad f(z) = \sum_{i=1}^{m} \alpha_i \frac{1}{\sigma_i} \phi\left(\frac{z - \mu_i}{\sigma_i}\right) \equiv \beta_1 f_1(z) + \beta_2 f_2(z)$$

where $\phi$ is the standard normal density, and the parameters are $\theta = (\alpha_1, \ldots, \alpha_m, \mu_1, \ldots, \mu_m, \sigma_1, \ldots \sigma_m)$, such that $0 < \epsilon < \sigma_i < \infty$, $\alpha_i > 0$ and $\sum_{i=1}^{m} \alpha_i = 1$. The first component density $f_1(z) = \phi((z - \mu_1)/\sigma_1)/\sigma_1$ represents the background distribution, while the second component density

$$(2) \qquad f_2(z) = \frac{1}{\beta_2} \sum_{i=2}^{m} \alpha_i \frac{1}{\sigma_i} \phi\left(\frac{z - \mu_i}{\sigma_i}\right)$$

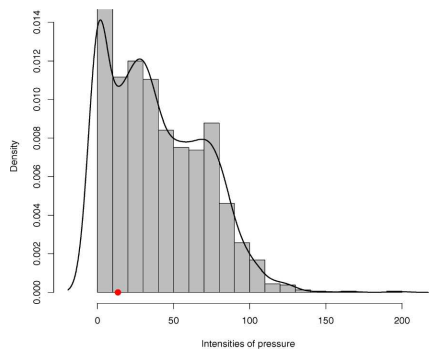

FIG 4. *Segmentation by analyzing the histogram and density plot of a data frame. A simple threshold is the red point which corresponds to the first deepest valley point between the first two consecutive major peaks in the density curve.*

is the signal distribution from the sitting region, which can be approximated well by a finite mixture of normal distributions. Here $\beta_1 = \alpha_1, \beta_2 = \sum_{i=2}^{m} \alpha_i$. See Figure 4 for the histogram and density plot of the data frame shown in the first subplot of Figure 2 (subject 1). Our analyses of data from NMES experiments showed that a mixture of two or three component normal distributions fitted our data quite well. In fact, the optimal estimate of $T$ developed below is fairly robust even if $f_2$ departs from a mixture of normal distributions slightly.

The following is our segmentation algorithm for determining the threshold $T$:

1). For a set of reasonable values of $m$, compute the estimates of $\theta$ by the *Expectation-Maximization* (EM) algorithm proposed by [6] for each given $m$. It is important to start from good initial values of $\theta$ in running an EM algorithm. We recommend choosing the initial values based on an under-smoothed histogram and summary statistics.
2). Estimate the final $m$ based on the Bayesian information criterion [13]. (In our NMES study, the final $m$ for all treatment subjects ended up to be 2 or 3.)
3). Select $T$ that is the solution of

$$(3) \qquad \alpha_1 \frac{1}{\sigma_1} \phi\left(\frac{T - \mu_1}{\sigma_1}\right) = \sum_{i=2}^{m} \alpha_i \frac{1}{\sigma_2} \phi\left(\frac{T - \mu_i}{\sigma_i}\right)$$

where $\alpha_i$, $\mu_i$ and $\sigma_i$ are the EM estimates from Step 1 and $m$ from Step 2.

That $T$ defined in (3) is optimal follows immediately from the following simple lemma by defining $f_1$ and $f_2$ as those in (2).

**Lemma 3.1.** *In a two-class classification problem assume that the probability density functions of the two populations are $f_i(z) \in C^2$, each has a mean $\mu_i$, for $i = 1, 2$, such that $\mu_1 < \mu_2 < \infty$. Furthermore, assume that the prior probability of population $i$ is $\beta_i > 0$, such that $\beta_1 + \beta_2 = 1$. Then the optimal threshold $T$ that minimizes the overall probability of misclassification (PMC) of a "simple" classification rule, such that an observation is classified into class one if $x \leq T$ and class two if $x > T$, is one such that $\beta_1 f_1(T) = \beta_2 f_2(T)$. Further, if $f_i$'s are normal densities with finite nonsingular variances, then the simple classification rule with $T$ defined above is also one that minimizes the overall PMC among all two-class classification rules.*

The proof of this lemma is straightforward. For example, for the general $f_i \in C^2$ case in which the classification rule is simple, one can simply write down the overall

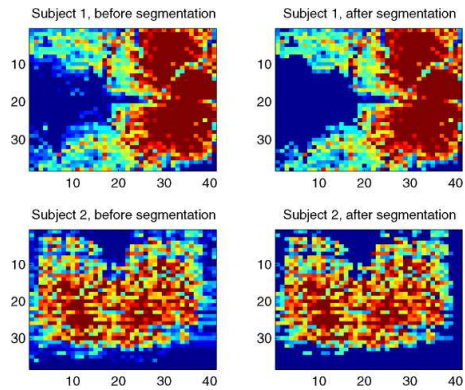

Fig 5. *Examples of comparison of images after data segmentation using optimal thresholds. The upper two subplots are for Subject 1 and the lower two subplots are for Subject 2. Note that the sitting regions in the data frames are segmented and the background noises are removed.*

PMC for a fixed $T$ and then differentiate the PMC with respect to $T$. The result for the mixture of normal densities can be obtain by tracing the equivalence between minimizing the overall PMC and maximizing the "posterior" probability, as shown in, for example, Result 11.5 in [12].

Figure 5 shows two examples of our data segmentation methods. The optimal thresholds in the data frames of subject 1 and subject 2 are 12.7 and 14.3, respectively. The sitting regions in the data frames are clearly segmented and the background noise is removed.

### 3.2. Data Registration

In the NMES study the image object was the anatomical seating contact area of the body, specifically the buttock and thigh region. The current experimental protocol entailed obtaining several data sets from each subject during their participation in the experiment. Since a subject may not sit at the same relative position on the sensor mat or with the same posture as previously, or the image target regions may differ from one session to another, some images from different sessions were not aligned. Recall Figure 2 where misalignment for some subjects is more obvious than the others. For example, the image in the second row for the fourth subject has been rotated 90 degrees in the last session. The image for subject 6 has non-overlapping areas between two images. Non-overlapping regions will be chopped out or trimmed during final analysis.

Since the subjects were not restrained in any way during the assessment it was also possible for some change in seating orientation to occur from one assessment condition to another during the same session. In order to determine any changes due to the effect of using NMES we first had to ensure that any changes due solely to seating orientation were fully compensated. This was achieved by spatial registration.

In the middle segment of each session as shown in Figure 1, a left/right alternating stimulation is given to a subject. To compare middle segments from two sessions a temporal registration is necessary to avoid artificial differences caused by stimulation cycle phase obscuring true image differences due to treatment.

*3.2.1. Spatial registration scheme: SRLP*

Generally, registration can be done by a geometrical transformation, which is a mapping of points from the space $\mathcal{A}$ of one view to the space $\mathcal{B}$ of a second view. The transformation **T** applied to a point in $\mathcal{A}$ represented by the column vector $\mathbf{a} = (a_i, a_j)^T$ produces a transformed point $\mathbf{a'} = (a'_i, a'_j)^T = \mathbf{a'} = \mathbf{T(a)}$. If the point $\mathbf{b} = (b_i, b_j)^T \in \mathcal{B}$ corresponds to **a**, then any nonzero displacement $\mathbf{T(a)} - \mathbf{b}$ is a registration error. Fortunately, the images within one segment in one session, and between different data sets in one session, do not appear to need spatial registration. Thus, we only need to spatially register images from different sessions. So, the first stable image of the first movie in each session can be used as a reference to register or align movies from different sessions, before we compute difference images or movies for statistical analysis of clinical relevance.

For spatial registration of data, a key is to choose appropriate landmarks. In our analysis of data from the NMES experiments, a natural landmark is the "midline" of the seating contact area for each patient. The midline and an obvious "end" point in each image will be used as our landmarks for registration leads to a midline-to-midline and endpoint-to-endpoint alignment. A scale change of images is not expected unless a subject has a significant change in body weight between two sessions. Thus we propose the following SRLP algorithm.

**Algorithm 3.1.** Automatic Spatial Registration by a line and a point (SRLP)

**1.** Determine the midpoints for each image,

$$\mathrm{midpt} = \frac{\mathrm{rowcount}}{2} + \frac{(c_1 - c_2)}{2}$$

where $c_1 =$ the number of non-zero values from the lower half image, $c_2 =$ the number of non-zero values from the upper half image, and rowcount is the total number of non-zero values in each column of the image.

**2.** Determine the midline. The midline is the regression line estimated by fitting a simple regression to the midpoints.

**3.** Perform a rigid transformation based on the midline, by rotation and translation through matrix **R**

$$\begin{bmatrix} a'_i \\ a'_j \\ 1 \end{bmatrix} = \mathbf{R} \begin{bmatrix} a_i \\ a_j \\ 1 \end{bmatrix} = \begin{bmatrix} \cos\theta & -\sin\theta & u \\ \sin\theta & \cos\theta & v \\ 0 & 0 & 1 \end{bmatrix} \begin{bmatrix} a_i \\ a_j \\ 1 \end{bmatrix}$$

where $\tan\theta$ is the slope of the midline and $(u, v)$ is the last point of the fitted midline in the image that is to be transformed.

If the patient is sitting asymmetrically the two halves of the image will have an unequal number of non-zero pixel values. For example, if the patient is leaning toward the lower half of the image there will be more non-zero pixel values in the lower half than in the upper half of the image, i.e $c_1 > c_2$. A positive correction $(c_1 - c_2)/2$ to the rowcount/2 should then be applied so that the location of the midpoint value moves up. After computation of the corrected midpoints, the midline can readily be found through linear regression. In Figure 6, the upper graph displays the midline of a patient in one frame; the lower graph displays the images after spatial registration for the same subject.

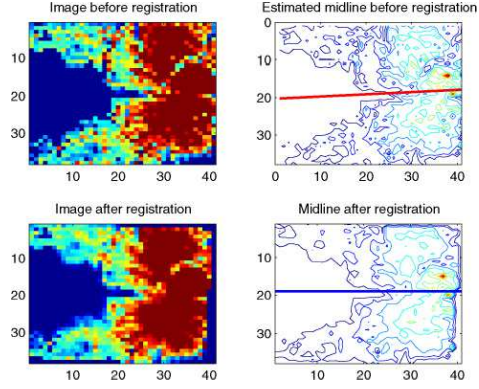

FIG 6. *An example of spatial registration by a line and a point (SRLP). The middle line is determined by a simple linear regression and rigid transformation is used in the registration.*

**Remark.** The idea of SRLP is simple but highly effective. It allows for self-registering any image based on its own midline found automatically by an algorithm. It can correct the bias and save the labor in determining the middle line manually. It is also a consistent algorithm in a statistical sense for a random landmark registration problem as shown in Theorem A.1 in the Appendix.

*3.2.2. A temporal registration scheme: ICR*

As part of the assessment protocol for this study electrical stimulation of the gluteal muscles was applied to produce dynamic weight-shifting from side to side. Temporal registration is required to align stimulation periods (on-off times) for all data sets collected for one subject under the same assessment conditions.

If the intensities in images A and B are linearly related, then the correlation coefficient is the ideal similarity measure. Few registration applications will precisely conform to this linear requirement, but many intra-modality applications, such as aligning on-off signals for two simulation sessions in our case, come sufficiently close for this to be an effective measure.

**Algorithm 3.2.** Intensity-based correlation registration (ICR)

**1.** Discard the first $m_0$ unstable data frames from each of the sub-data sets with the NMES stimulation (Here we choose $m_0 = 10$).

**2.** For the remaining images $A_1, \ldots, A_n$ and $B_1, \ldots, B_n$ from the middle segments of two on-off stimulation sessions, compute the correlation coefficient $cor_{ij}(AB)$ of $A_i$ and $B_{i+j}$ for $i = 1, \ldots, n-j$ and $j = 0, \ldots, n-1$. Let

$$\text{CorAvg}_j = \frac{1}{n-j} \sum_i cor_{ij}(AB).$$

Find $j_0$ such that

$$\text{CorAvg}_{j_0} = \max_j(\text{CorAvg}_j).$$

**3.** Align images $A_i$ with $B_{i+j_0}$.

## 4. Statistical smoothing mapping

Our primary questions of interest in the NMES study are: 1) Does the long-term gluteal NMES improve intrinsic characteristics of the paralyzed muscles? 2) Can we identify the areas in which interface pressure has significantly improved?

In statistical methods of brain imaging (e.g. MRI), one of the most common analysis approaches currently in use, called *statistical parametric mapping* (SPM) [8, 9], analyzes each voxel's change independently of the others and builds a map of statistic values for each voxel. The significance of each voxel can be ascertained statistically with a Student's t-test, an F-test, a correlation coefficient, etc. SPM is widely used to identify functionally specialized brain regions and is the most prevalent approach to characterizing functional anatomy and disease-related changes. The success of SPM is due largely to the simplicity of the idea. Namely, one analyzes each and every voxel using any standard (univariate) statistical parametric test. The resulting statistical parameters are assembled into an image – the SPM.

Motivated by the SPM, we propose a statistical smoothing mapping (SSM) procedure based on multivariate smoothing, to allow for more flexible modeling than parametric models. Since we are comparing many voxel values *simultaneously* across the entire image, the multiplicity of these tests must be adjusted to overcome an overall false-positive error rate. Our significance threshold for deciding which voxel is significantly different (between two sessions) will be chosen with a Benjamini and Yekutieli false discovery rates (BH-FDR) controlling procedure [1] that accounts for the multiplicity of tests. Then an FDR map can be built to provide the significance of voxels. Those with p-values less than the BH critical value are the points or areas for which stimulation has had a significant effect (difference) in terms of measurements.

Let $\tilde{x} = (x_1, x_2)$ denote a pixel of a data frame. Then $r_{\tilde{x},C}$, $r_{\tilde{x},T}$ denote the intensities of the (registered) images before treatment and after treatment. We propose the following statistical smoothing mapping algorithm.

**Algorithm 4.1.** Statistical Smoothing Mapping (SSM)

1. Compute the difference map, $y_{\tilde{x}} = r_{\tilde{x},T} - r_{\tilde{x},C}$ which is the pixel-by-pixel subtraction before treatment and after treatment. Then pad the same values of $y_{\tilde{x}}$ at the edge of sitting regions into a small rim of the background region to overcome the possible edge effects of smoothing techniques in the next step.
2. Smooth padded $y_{\tilde{x}}$ by multivariate local polynomial regression.
3. Compute the "t-type" statistic $T_{\tilde{x}}$ (defined below) and p-value for each pixel. Then chop off the "padded t"-values outside the sitting region.
4. Compute adjusted p-values using the BH-FDR controlling procedure. Generate an FDR map/movie based on the adjusted p-values.

### *4.1. T-type tests*

In the SSM algorithm, we consider the following nonparametric regression model to smooth $y_{\tilde{x}}$,

$$Y_i = m(\tilde{X}_i) + \varepsilon_i$$

where $\tilde{X} = (X_1, X_2)$ is a two-dimensional predictor which denotes the coordinates of an image; the response variable $Y_i$ is the corresponding intensity at $\tilde{x}$; $m(\cdot, \cdot)$ is an unknown smooth function and $\varepsilon_i$ is an error term, representing random errors in the observations and variability from sources not included in the $\tilde{X}_i$. A smooth function $m$ can be approximated in a neighborhood of a point $\tilde{x} = (x_1, x_2)$ by a local polynomial. Here we consider a local quadratic approximation:

$$m(\tilde{u}|\tilde{a}) \approx a_0 + a_1(u_1 - x_1) + a_2(u_2 - x_2)$$
$$+ \frac{a_3}{2}(u_1 - x_1)^2 + \frac{a_4}{2}(u_2 - x_2)^2 + a_5(u_1 - x_1)(u_2 - x_2)$$

where $\tilde{u} = (u_1, u_2)$, $\tilde{a} = (a_1, a_2, a_3, a_4, a_5, a_6)$. The coefficient vector $\tilde{a}$ can be estimated by minimizing a locally weighted sum of squares:

$$\sum_{i=1}^{n} w_i(\tilde{x})(Y_i - m(\tilde{x}_i|\tilde{a}))^2$$

where $w_i(\tilde{x})$ is a spherically symmetric weight function that gives an observation $\tilde{x}_i$ the weight $w_i = W(||\tilde{x}_i - \tilde{x}||/h)$. The local regression estimate of $m(\tilde{x})$ is defined as $\hat{m}(\tilde{x}) = \hat{a}_0$. (See [17] about computational aspects and bandwidth selections in details.)

Our hypotheses at $\tilde{x}$ here are: $H_0 : m(\tilde{x}) = 0$ vs. $H_1 : m(\tilde{x}) > 0$. Since $\hat{m}(\tilde{x})$ can be written as a linear combination of the response variables,

(4) $$\hat{m}(\tilde{x}) = \sum_{i=1}^{n} p_i(\tilde{x})Y_i,$$

where $p(\tilde{x})^T = (p_1(x), \ldots, p_n(x))$ is the rows of the *hat matrix* specified by the quadratic approximation, the estimated standard deviation of the local estimate $\hat{m}$ is $\hat{S}(\tilde{x}) = \hat{\sigma}||p(\tilde{x})||$. A proper test statistic is then the "t-type" statistics:

(5) $$T(\tilde{x}) = T_{\tilde{x}} = \frac{\hat{m}(\tilde{x})}{\hat{S}(\tilde{x})}.$$

The null hypothesis $H_0$ would be rejected at $\tilde{x}$ if $T(\tilde{\mathbf{x}}_i) > t_{1-\alpha}(\delta_1^2/\delta_2)$ with a given significance level $\alpha$, where the degrees of freedom $\delta_1$ and $\delta_2$ can be obtained by two-moment chi-square approximations [5, 16, 17], if we used a pointwise test.

Since the above test statistic $T$ is a weighted average of $y$ values in a neighborhood of $x$, $T(\tilde{x})$ and $T(\tilde{x}')$ are often correlated if $\tilde{x}$ and $\tilde{x}'$ are not far away.

### 4.2. Multiple testing problem

For multiple comparisons, an overall error must be controlled to overcome the *multiplicity* problem that occurs with simultaneously testing many hypotheses, $H_0 : m(\tilde{x}) = 0$ for all $\tilde{x}$ in the sitting region. The family-wise error rate (FWER) and false discovery rate (FDR) are two typical overall error rates. The simplest multiple comparison procedure that controls the FWER is the Bonferroni procedure. However, the Bonferroni procedure is too conservative when the number of hypotheses is very large. It is important to note that the conservativeness of the Bonferroni procedure comes from two sources: (1) the Bonferroni procedure was based on a very conservative *upper bound* for the FWER. (2) FWER is a more stringent error than FDR. To overcome (1), there are sharper upper bounds for FWER developed for finite $m$ cases (see [11] and reference therein); there are also exact and accurate approximations to FWER by tube formulas such as those shown in [14, 15, 16, 18]. As to when to use FWER or FDR, see Section 4.2 of [19]. In this paper, we choose FDR.

The step-up procedure for strong control of the false discovery rate introduced by [1] can be easily implemented, even for very large data sets. We call this BH-FDR procedure. Returning to the NMES study, it is observed that the approximate T statistics of the multiple tests are dependent as they are from the estimated regression function. [2] showed that the BH-FDR procedure is valid under "positive regression dependency on subsets" (PRDS). They also proposed a simple conservative modification of the procedure which controls the false discovery rate for arbitrary dependence structures by letting $c_m = \sum_{i=1}^{m} 1/i$. Note that $\sum_{i=1}^{m} 1/i \approx \ln m + \gamma$ where $\gamma$ is the Euler's constant. For a large number $m$ of hypotheses, the penalty in this conservative procedure is about $\log m$, (as compared to the BH-FDR procedure) which can be still too large and can be more conservative than the tube methods or random field methods by [15, 16]. Rather than using this conservative procedure with a factor $\ln m + \gamma$, we prove that the joint distribution of the T-type test statistics in multivariate local regression is PRDS on the subset of test statistics corresponding to true null hypotheses, and thereby the BH-FDR procedure is still valid – See Appendix B.

## 5. LASR – A new data mining procedure

### 5.1. LASR

Combining the techniques we developed, we present a complete data-mining scheme, the LASR procedure for analyzing a large sequence of spatial-temporal data sets.

Figure 7 displays the flow chart of our LASR procedure. Step 1: Perform *Segmentation* for all images by the EM algorithm and by computing the optimal threshold. Segment the spatial regions of interest from the background in each data frame and then remove background noise and outliers from the data sets. Segmentation is needed only for one image frame per movie. Step 2: *Spatially register* all (segmented) images via our self-registration scheme SRLP. This step is done automatically for all images so that all registered images have the middle line placed horizontally in the middle of each image and the end point at the same location. If both movies are static movies, go to step 3; if both are dynamic movies, *temporally register* the spatially-registered movies. The temporal registration is based on the ICR algorithm to maximize the correlations between images from two candidate movies, frame-by-frame so that the left side that is stimulated in one movie is compared with the left-side stimulated image in another movie (See movies at *stat.case.edu/lasr/*). Step 3: *Create difference images and movies* by taking differences pixel-by-pixel (and frame-by-frame) between two sessions that are potentially of clinical interest. Step 4: *Smooth* the difference images by the bivariate nonparametric local regression as in Section 4.1. Step 5: *Create T image maps and movies*. Generate T images by computing a t-type test statistic at each pixel. Step 6: *Compute FDR-controlled P maps and movies*. Based on the T images and movies, we can compute (pointwise) $p$-values at all pixels. The BH-FDR method is applied to adjust the p-values to account for the multiplicity from simultaneously testing for differences at all pixels. If a p-value $p$ at $x$ is less than the critical value derived from a 0.05 FDR-controlled procedure, change the pixel value to $1 - p$; if $p$ is greater than the FDR cut-off value, the pixel value is set to zero. These resulting FDR-controlled P maps or movies show which areas are the elevated areas or the area that improve interface pressures (implying improved tissue health).

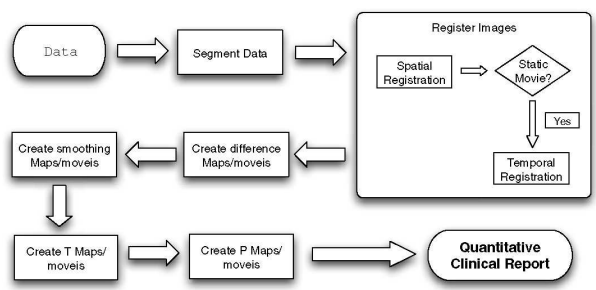

FIG 7. *LASR procedure flow chart.*

The LASR output map gives a graphical representation of statistically significant pressure changes across the entire mapped region. It helps us to decide if the NMES is effective at a particular region, with an FDR no more than 0.05. The algorithm is applied frame-by-frame to aligned pressure data sets. LASR maps can thus be viewed as single frame "snapshots", suitable for comparison of static seating postures, or as videos for comparison of dynamically changing pressures.

### 5.2. Statistical results

In this subsection we present three typical analysis results: one for a control subject who did not receive any NMES, one for a treatment subject with static mappings and another for a treatment subject with dynamic mappings. The results for other five subjects support the conclusion we can draw from these three analysis results.

*Control case: Subject A*. Seating pressure assessments were obtained for subject A at an interval of three months, during which time no NMES was used. From the upper subplot of Figure 8, it is noted that some spatial misalignment is apparent between baseline and repeat assessment. After applying the LASR algorithm, it could be seen from the blank P-map that there was no significant differences in interface pressure distributions obtained at a three month interval for an individual who was not receiving NMES.

*NMES users: Subject B in static model, subject C in dynamic model*. Seating pressure assessments were obtained for subject B and subject C at an interval of six months, during which time NMES was used regularly. After applying the LASR algorithm to assess changes between baseline and post-treatment interface pressure data sets it could be seen that pressures were reduced bilaterally over time. Figure 9 (b) shows long-term changes for subject B in static mode seated pressure distribution. The left sacro-ischial region was more extensively affected than the right side. Figure 9 (c) shows long-term changes for subject C in dynamic mode seated pressure distribution. The left and right sacro-ischial regions were equally affected. Relevant LASR movies can be viewed at *stat.case.edu/lasr/*.

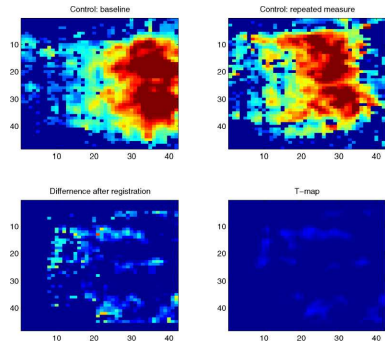

Fig 8. *LASR analysis results for control data.*

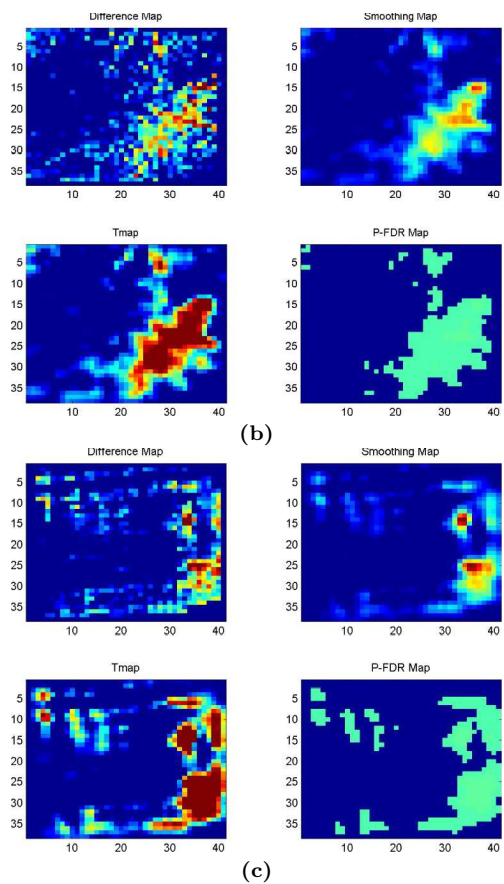

FIG 9. *LASR analysis for pressure mappings to identify the regions of the pressure reduction.* **(b)** *(upper 4 pictures): Subject B in static mode;* **(c)** *(lower 4 pictures): Subject C in dynamic mode.*

## 6. Discussion

The development of the multi-stage statistical LASR algorithm allows both clinicians and researchers to derive useful, objective information from pressure maps, such as the location of significant pressure changes or the relative efficacy of pressure relief procedures. Furthermore, spatial registration allows global analysis of pre- and post-intervention differences without any subjective bias in selecting areas of interest. In the specific study of the effects of NMES it was found that subjects who received a gluteal stimulation system showed statistically significant changes in ischial region pressure over time, when baseline/post-treatment comparisons were made. The region of significant change was not symmetrical in all cases which re-flects the asymmetric nature both of gluteal muscle recruitment area and contractile responses.

The last two decades have seen remarkable developments in imaging technology. Medical images are increasingly widely used in health care and biomedical research; a wide range of imaging modalities is now available. The clinical significance of medical imaging in the diagnosis and treatment of diseases is overwhelming. The commonly used medical imaging modalities capable of producing multidimensional images for clinical applications are: X-ray Computed Tomography (X-ray CT), Magnetic Resonance Imaging (MRI), Single Photon Emission Computed Tomography (SPECT), Positron Emission Tomography (PET) and Ultrasound (US).

It should be noted that these modern imaging methods involve sophisticated instrumentation and equipment which employ high-speed electronics and computers for data collection. Spatial-temporal (image) data occur in a broad range of medical applications. It is now common for patients to be imaged multiple times, either by repeated imaging with a single modality, or by imaging with different modalities. It is also common for patients to be imaged *dynamically*, that is, to have sequences of images acquired, often at many frames per second. The ever increasing amount of image data acquired makes it more and more desirable to relate more than one statistical tool to assist in extracting relevant clinical information.

Application of the LASR algorithm enhances data extraction and acquires statistical inferences from complex spatial-temporal data sets, as shown in the NMES study. Thus the LASR analytical methodology has the potential to be applied to other imaging modalities or to other imaging targets in which natural landmarks may be different from the midline and an end point. Other potential clinical applications include images of soft tissues, which may not include bony landmarks. Applications could include situations where an imaged object may change dimensions and/or orientation over time.

## Appendix A: Registration error of SRLP

Define the overall *registration error* (RE) of a transformation $\mathbf{T}$ to be

$$\text{RE} = \frac{1}{n^2} \sum_{i=1}^{n} \sum_{j=1}^{n} ||\mathbf{T}(\mathbf{a}_{ij}) - \mathbf{b}_{ij}||^2 \tag{6}$$

where $\mathbf{a}_{ij}$ and $\mathbf{b}_{ij}$ ($i, j = 1, \ldots, n$) are the corresponding points (*i.e.* pixel coordinates) in spaces $\mathcal{A}$ and $\mathcal{B}$, respectively.

**Theorem A.1.** *Assume that the intensity values are bounded and are confined in a bounded domain. Then the SRLP is consistent in terms of RE as the number of pixels increases.*

*Proof.* After SRLP registration, the $\mathbf{a}'_{ij} = (a'_i, a'_j)$ has the representation

$$a'_i = a_i \cos \hat{\theta} - a_j \sin \hat{\theta} + u, \qquad a'_j = a_i \sin \hat{\theta} + a_j \cos \hat{\theta} + \hat{v}$$

where $\tan \hat{\theta} = \hat{\beta}_0$, $\hat{v} = \tan \hat{\theta} u + \hat{\beta}_1$, and $\hat{\beta}_0$, $\hat{\beta}_1$ are the estimates of the slope and intercept of the midline. Notice that $u$ is not an estimated value because the horizontal axis of the last point in the fitted midline keeps immovable.

A perfect registration will make the transformed point equal to $\mathbf{b}_{ij}$, so

$$b'_i = a_i \cos \theta - a_j \sin \theta + u, \qquad b'_j = a_i \sin \theta + a_j \cos \theta + v$$

Then the registration error of SRLP is equal to

$$\text{RE} = \frac{1}{n^2} \sum_{i=1}^{n} \sum_{j=1}^{n} \Big\{ [a_i(\cos\hat{\theta} - \cos\theta) - a_j(\sin\hat{\theta} - \sin\theta)]^2 \\ + [a_i(\sin\hat{\theta} - \sin\theta) - a_j(\cos\hat{\theta} - \cos\theta) + (\hat{v} - v)]^2 \Big\}.$$

Note that $\hat{\beta}_0$, $\hat{\beta}_1$ are consistent estimators in the midline regression, and $\theta = g(\beta_1, \beta_2)$, $v = h(\beta_1, \beta_2)$ where $g$ and $h$ are continuous functions. Hence, $\hat{\theta} = g(\hat{\beta}_0, \hat{\beta}_1)$ and $\hat{v} = h(\hat{\beta}_0, \hat{\beta}_1)$ are also consistent by *Slutsky's theorem*. Then by the boundedness of intensities and $a_i$, $a_j$, it is easy to see that $\text{RE} \to 0$ in probability as $n \to \infty$ or the number of pixels tends to infinity. □

### Appendix B: PRDS of test statistics in multivariate local regression

Recall that a set D is called increasing if $x \in D$ and $y \geq x$, implies that $y \in D$ as well. The following property is called *positive regression dependency on each one from a subset* $I_0$, or PRDS on $I_0$ [2].

**Property B.1** (PRDS). For any increasing set $D$, and for each $i \in I_0$, $P(X \in D|X_i = x)$ is nondecreasing in $x$.

**Proposition B.1** (PRDS of test statistics in multivariate local regression). *Consider a vector of test statistics* $\mathbf{T} = (T_1, T_2, \ldots, T_m)^T$. *Each* $T_i$ *tests the hypothesis* $m(\tilde{x}_i) = 0$ *against the alternative* $m(\tilde{x}_i) > 0$ *for* $i = 1, \ldots, m$, *where* $T_i$ *is defined by* (5) *with the nonnegative weights* $p(\tilde{x}_i)$ *in* (4). *The distribution of* $\mathbf{T}$ *is PRDS over* $I_0$, *the set of true null hypotheses.*

*Proof.* Let $\mathbf{U} = (U_1, \ldots, U_m)^T$ where $U_i = \hat{m}(\tilde{x}_i)/||p(\tilde{x}_i)||$. We first verify that $\mathbf{U}$ is PRDS on a subset $I_0$. By (4), for any $i \neq j$,

$$\text{cov}(U_i, U_j) = \frac{\text{cov}\big(\sum_{t=1}^{n} p_t(\tilde{x}_i) Y_t, \sum_{k=1}^{n} p_k(\tilde{x}_j) Y_k\big)}{||p(\tilde{x}_i)|| \cdot ||p(\tilde{x}_j)||}$$
$$= \frac{\sigma^2 \sum_{t=1}^{n} p_t(\tilde{x}_i) p_t(\tilde{x}_j)}{||p(\tilde{x}_i)|| \cdot ||p(\tilde{x}_j)||} > 0$$

Under the normality assumption of errors, $\mathbf{U}$ follows a multivariate normal distribution with the covariance matrix having positive elements. Then $\mathbf{U}$ is PRDS on a subset $I_0$ because the conditional distribution $\mathbf{U}_{(i)}$ given $U_i = u_i$ increases stochastically as $u_i$ increases (where $\mathbf{U}_{(i)}$ denotes the remaining $m-1$ test statistics except $U_i$).

Since $\hat{\sigma}^2$ approximately follows $\chi^2$ distribution, let $V = 1/\hat{\sigma}$. Then for $j = 1, \ldots, m$ the components of $\mathbf{T}$, $T_j = U_j V$ are strictly increasing continuous functions of the coordinates $U_j$ and of $V$. Therefore, $\mathbf{U}$ is PRDS on $I_0$ by applying Lemma 3.1 of [2]. □

### Acknowledgments

Research is supported in part by a grant from the DMS in NSF.